\documentclass[11pt]{amsart}
\newtheorem{theorem}{Theorem}[section]
\newtheorem{corollary}[theorem]{Corollary}
\newtheorem{lemma}[theorem]{Lemma}
\newtheorem{proposition}[theorem]{Proposition}

\theoremstyle{definition}

\theoremstyle{remark}

\numberwithin{equation}{section}

\newcommand{\C}{\mathbb C}

\begin{document}

\title{ Extension of CR maps of positive codimension}

\author{Sergey Pinchuk and Alexandre Sukhov}

\subjclass[2000]{32H02, 53C15}

\date{\number\year-\number\month-\number\day}

\begin{abstract}
We study the holomorphic extendability of smooth CR maps between real analytic strictly pseudoconvex hypersurfaces in complex affine spaces of different dimensions. 

\end{abstract}

\maketitle

\section{Introduction}
This paper concerns the following long-standing conjecture: let $f: M \longrightarrow M'$ be a smooth CR map 
between two real analytic strictly pseudoconvex hypersurfaces in the complex affine spaces $\C^n$ and $\C^N$ respectively with $1 < n \leq N$.
Then $f$ extends holomorphically to a neighborhood of $M$. At present  the strongest result is due to Forstneric \cite{Fors} who proved that
 $f$ extends to a neighborhood of an open dense subset of $M$. 

Here  we prove the following

\begin{theorem}
Let $M \subset \C^n$, $M' \subset \C^N$ be $C^{\omega}$ strictly pseudoconvex hypersurfaces and $f: M \longrightarrow M'$ be a $C^{\infty}$ 
CR map. If $2 \leq n \leq N < 2n$ then $f \in {\mathcal O}(M)$.
\end{theorem}

This gives a complete solution to the above problem in the case where the ``codimension'' $N -n$ of the map $f$ is smaller than $n$. 

{\it Aknowledgements.} This work was initiated when the first author visited the Universities Aix-Marseille-I and Lille-I in May 2004 and 
completed when the second author visited the University of Illinois during the Spring semester 2005. They thank  these institutions.

\section{ Notations and preliminaries}

Denote by $z = (z_1,...,z_n) \in \C^n$ and $z' = (z'_1,...,z'_N) \in \C^N$ the standard coordinates in $\C^n$ and $\C^N$ 
respectively. Without loss of generality we may assume that $0 \in M$, $0' \in M'$ and $f(0) = 0'$. It is enough to prove that $f$ extends
 holomorphically to a neighborhood of the origin. 

Consider sufficiently small connected neighborhoods ${\mathcal U}$ and
${\mathcal U'}$ of $0$ and $0'$ respectively. 
Let  $\rho(z) \equiv \rho(z,\overline z) \in
C^{\omega}({\mathcal U})$ and $\rho'(z') \equiv \rho'(z',\overline z') \in C^{\omega}({\mathcal U'})$ be strictly plurisubharmonic defining functions of $M$ and $M'$ respectively.
We will denote by $\rho(z,\overline w)$, $\rho'(z',\overline w')$
their complexifications.
If $\omega = \overline w$, $\omega' = \overline w'$, then 
$\rho(z,\omega) \in {\mathcal O}({\mathcal U} \times {\mathcal U})$,
$\rho'(z',\omega') \in {\mathcal O}({\mathcal U'} 
\times {\mathcal U'})$. 
 
For $w \in {\mathcal U}$ denote by $Q_w:= \{ z \in {\mathcal U}:
\rho(z,\overline w) = 0 \}$ 
the Segre variety of $w$. The Segre variety $Q'_{w'}$ is defined similarly
for $w' \in {\mathcal U}'$. Consider also the one-sided neighborhoods 

\begin{eqnarray*}
& &{\mathcal U}^{+}:= \{ z \in {\mathcal U}: \rho(z) > 0 \}, {\mathcal U}^{-}:= \{ z \in {\mathcal U}: \rho(z) < 0 \},\\
& &{\mathcal U'}^{+}:= \{ z' \in {\mathcal U'}: \rho'(z') > 0 \}, {\mathcal U'}^{-}:= \{ z' \in {\mathcal U'}: \rho'(z') < 0 \}
\end{eqnarray*}
 Then $f$ extends holomorphically to ${\mathcal U}^-$, and we may assume that  $f({\mathcal
   U}^-) \subset {\mathcal U'}^-$, 
$f \in C^{\infty}({\mathcal U}^- \cup M)$. Furthermore,
by Forstneric \cite{Fors} there exists an open dense subset $\Sigma
   \subset M \cap U$ such that 
$f \in {\mathcal O}({\mathcal U}^- \cup \Sigma)$. 

If $a \in \Sigma$ and $f$ is holomorphic on a neighborhood $V$ of $a$,
then $\rho'(f(z)) \in C^{\omega}(V)$ and 
$\rho'(f(z)) = \alpha(z)\rho(z)$ for 
$\alpha(z) \in C^{\omega}(V)$. After the complexification we have 
$\rho'(f(z),\overline{f(w)}) = \alpha(z,\overline w) \rho(z,\overline
w)$. 
This implies that for $w$ close enough to $a \in \Sigma$ we have

\begin{eqnarray}
\label{1}
f(Q_w \cap V) \subset Q'_{f(w)}
\end{eqnarray}

Thus, if $f$ extends holomorphically across $M$, then the graph of the
extended map $f$ over ${\mathcal U}^+$ 
must be contained in the set 

\begin{eqnarray}
\label{2}
F:= \{ (w,w') \in {\mathcal U}^+ \times {\mathcal U}': f(Q_w \cap {\mathcal U}^-) \subset Q'_{w'} \}
\end{eqnarray}
(Notice that since $M$ is strictly pseudoconvex and ${\mathcal U}$ is a
small neighborhood of the origin, then 
$Q_w \cap {\mathcal U}^-$ is connected, see \cite{DiPi1,DiPi2}).

The set $F$ has already been used by Forstneric in \cite{Fors} and our proof of Theorem is based on the result of \cite{Fors} and further careful study of $F$.

If $d(z,M)$ denotes the euclidean distance from $z \in {\mathcal U}$ to $M$, then it is wellknown that if the map $f$ is {\it non-constant} then for $z \in {\mathcal U}^-$ 

\begin{eqnarray}
\label{3}
d(f(z),M') \sim d(z,M)
\end{eqnarray}
which here and later means that there exists a constant $c > 0$ such that

\begin{eqnarray}
\label{3'}
\frac{1}{c}d(z,M) \leq d(f(z),M') \leq c d(z,M)
\end{eqnarray}
for all $z \in {\mathcal U}^-$. The left part of (\ref{3'}) is the consequence of the Hopf lemma while the right part follows 
from the assumption $f \in C^{\infty}(M)$. Another wellknown fact is that in this case the differential $df$ has maximal rank near $M$, 
i.e. we may assume that $f: {\mathcal U}^- \longrightarrow {\mathcal U}'$ is an embedding.

Consider a $C^{\infty}$ extension of $f$ to ${\mathcal U}$ which we
denote  by $\tilde f$. 
We may assume that $\tilde f: {\mathcal U} \longrightarrow {\mathcal U}'$ is a proper embedding and thus
$S':= \tilde f({\mathcal U}) \subset {\mathcal U}'$ is a closed $C^{\infty}$ manifold extending $f({\mathcal U}^-)$. 

\begin{lemma}
\label{lemma1.1}
Let $\rho$ be a strictly plurisubharmonic $C^\infty$ function on ${\mathcal U}$. Then (after shrinking ${\mathcal U}$) 
\begin{eqnarray}
\label{4}
\rho(z,\overline z) + \rho(w,\overline w) - \rho(z,\overline w) - \rho(w,\overline z) \sim \vert z - w \vert^2
\end{eqnarray}
for $z,w \in {\mathcal U}$.
\end{lemma}
\proof Let $\rho(z,\overline z) = \sum_{k,l}c_{k,l}z^k\overline z^l$ be the Taylor expansion of $\rho$ at $0$ in the multi-index
notation. Since $\rho$ is a real function, we have $c_{lk} = \overline c_{kl}$. Let also 
\begin{eqnarray*}
\phi(z,\overline w) := \sum c_{kl}(z^k\overline z^l + w^k \overline w^l - z^k \overline w^l - w^k \overline z^l) = \sum
c_{kl}(z^k - w^k)(\overline z^l - \overline w^l)
\end{eqnarray*}

Then
\begin{eqnarray*}
\rho(z,\overline z) + \rho(w,\overline w) - \rho(z,\overline w) - \rho(w,\overline z) = 
\phi(z,\overline w)=  L(z-w) + o(\vert z - w \vert^2)
\end{eqnarray*}
where $L(z-w)$ here denotes the Levi form of $\rho$ at $0$ which satisfies $L(z-w) \sim \vert z - w \vert^2$.

\begin{corollary}
\label{cor1}
In the situation of lemma \ref{lemma1.1} there exists a constant $c > 0$ such that for any $z \in {\mathcal U}^-$ 
we have the inclusion $Q_z \cap {\mathcal U} \subset {\mathcal U}^+$ and
\begin{eqnarray}
\label{5}
d(Q_z \cap {\mathcal U}, M) \geq c d(z,M)
\end{eqnarray}
\end{corollary}
\proof For any $w \in Q_z$ we have $\rho(z,\overline w) = \rho(w,\overline z) = 0$. By (\ref{4}) we have $\rho(z,\overline z) + \rho(w,\overline w) \geq 0$.
Since $d(z,M) \sim \vert \rho(z,\overline z) \vert$ and $\rho(z,\overline z) < 0$ this implies  $\rho(w,\overline w) > 0$ (i.e. $w \in {\mathcal U}^+$) and
$d(w,M) \geq c d(z,M)$.

\begin{corollary}
\label{cor2}
If that $z \in M \cap {\mathcal U}$ and $w \in Q_z \cap {\mathcal U}$, then 
\begin{eqnarray*}
d(w,M) \sim \vert w - z \vert^2
\end{eqnarray*}
\end{corollary}
This directly follows from (\ref{4}).

\begin{lemma}
\label{lemma1.4}
$F$ is an analytic set in ${\mathcal U}^+ \times {\mathcal U}'$ of dimension $\geq n$.
\end{lemma}

\proof We assume that $\frac{\partial \rho}{\partial z_1}(0,0) \neq 0$ and therefore for any $w \in {\mathcal U}$ 
the equation $\rho(z,\overline w) = 0$ of $Q_w$ is equivalent to $z_1 = h(z_2,...,z_n,\overline w)$, where $h$ is holomorphic in $z_2,...,z_n$
and antiholomorphic in $w$. Thus for $w \in {\mathcal U}^+$ the condition $f(Q_w \cap {\mathcal U}^-) \subset Q'_{w'}$ is equivalent to the condition that for every
$z=(z_1,...,z_n) \in Q_w \cap {\mathcal U}^-$
\begin{eqnarray}
\label{5}
\rho'(f(h(z_2,...,z_n,\overline w),z_2,...,z_n), \overline w') = 0
\end{eqnarray}
This is a system of (anti)holomorphic equations for $w,w'$. Since $F$
is obviously closed in 
${\mathcal U}^+ \times {\mathcal U}'$, it is an analytic set. If $w,z$ are close to 
a point of holomorphic extendability of $f$, then $\rho(z,\overline w) = 0$ implies $\rho'(f(z),\overline{f(w)}) = 0$ and thus $F$ contains a piece of the graph of the extension of $f$ and $\dim F \geq n$.

\section{Boundary behaviour of $F$}

Let $\overline F$ be the closure of $F$ in ${\mathcal U} \times {\mathcal U}'$ and $\pi$, $\pi'$ be the natural projections of ${\mathcal U} \times {\mathcal U}'$ to 
${\mathcal U}$ and ${\mathcal U}'$ respectively. 

First notice that if $(w,w') \in \overline F$ and $w \in M$, then $w' \in Q'_{f(w)}$. Indeed, let $(w,w') = \lim_{\nu \longrightarrow \infty} (w^{\nu},{w'}^{\nu})$,
$(w^{\nu},{w'}^{\nu}) \in F$, and $z^{\nu} \in Q_{w^{\nu}} \cap {\mathcal U}^-$. Consider a sequence $z^{\nu} \in Q_{w^{\nu}} \cap {\mathcal U}^-$ such that $z^{\nu} \longrightarrow 
w$ and so $f(z^{\nu}) \longrightarrow f(w)$. Then $f(z^\nu) \in
Q'_{{w'}^\nu}$ and 
${w'}^\nu \in Q'_{f(z^\nu)} \longrightarrow Q'_{f(w)}$ so that $w' \in Q'_{f(w)}$. This can
be reformulated as

\begin{eqnarray}
\label{6}
\overline F \cap (\{ w \} \times {\mathcal U}') \subset \{ w \}\times Q'_{f(w)}
\end{eqnarray}
for $w \in M \cap {\mathcal U}$.

We will now improve (\ref{6}). Differentiating (\ref{5}) with respect to $z_k$, $k=2,...,n$ we get

\begin{eqnarray*}
\sum_{j=1}^N \rho'_j(f(z),\overline w')\left ( \frac{\partial f_j}{\partial z_1}(z)\frac{\partial h}{\partial z_k}(\tilde z,\overline w) + 
\frac{\partial f_j}{\partial z_k}(z) \right ) = 0
\end{eqnarray*}
where $\tilde z = (z_2,...,z_n)$ and $\rho'_j:= \frac{\partial \rho'}{\partial z'_j}$. Since 
$$\frac{\partial h}{\partial z_k}(\tilde z,\overline w) = - \frac{\rho_k(z,\overline w)}{\rho_1(z,\overline w)}$$
for $z \in Q_w \cap {\mathcal U}$, this is equivalent to 

\begin{eqnarray}
\label{7}
\sum_{j=1}^N \rho'_j(f(z),\overline w') T_kf_j(z,\overline w) = 0, k=2,...,n
\end{eqnarray}
 where 
\begin{eqnarray}
\label{8}
T_kf_j(z,\overline w):= \rho_1(z,\overline w) \frac{\partial f_j}{\partial z_k} - \rho_k(z,\overline w) \frac{\partial f_j}{\partial z_1}(z)
\end{eqnarray}

In particular, if $w \in M$ we can take $z = w$ and (\ref{7}) becomes

\begin{eqnarray*}
\sum_{j=1}^N \rho'_j(f(w),\overline w') T_kf_j(w,\overline w) = 0, k=2,...,n
\end{eqnarray*}

Thus we proved

\begin{lemma}
\label{lemma2.1}
If $w \in M \cap {\mathcal U}$ then
\begin{eqnarray*}
\overline F \cap (\{ w \} \times {\mathcal U}') \subset \{ w \} \times \{ w' \in {\mathcal U}' \cap Q'_{f(w)}: 
\sum_{j=1}^N \rho'_j(f(w),\overline w') T_kf_j(w,\overline w) = 0, k=2,...,n \}
\end{eqnarray*}
\end{lemma}

Consider now a $C^{\infty}$ extension $\tilde f$ of $f$ to ${\mathcal U}$. Since $df$ has maximal rank at $0$ we may assume that it remains maximal in ${\mathcal U}$ and $\tilde f$ is a proper embedding of ${\mathcal U}$ to ${\mathcal U}'$. The image $S' = \tilde f({\mathcal U}) \subset {\mathcal U}'$ is a $C^\infty$ manifold of real dimension $2n$ which 
extends $f({\mathcal U}^-)$. 

\begin{lemma}
\label{lemma2.2}
For $(w,w') \in \overline F$ with $w \in M \cap {\mathcal U}$ 
\begin{eqnarray}
\label{9}
d(w',S') \sim \vert w' - f(w) \vert
\end{eqnarray}
\end{lemma}
\proof Choose the local coordinates near $ 0 \in \C^n$ and $0' \in \C^N$ such that 
\begin{eqnarray}
\label{10}
\rho(z) = 2x_1 + \vert z \vert^2 + o(\vert z \vert^2), \rho'(z') = 2x'_1 + \vert z' \vert^2 + o(\vert z'\vert^2)
\end{eqnarray}
\begin{eqnarray}
\label{11}
f_j(z) = z_j + o(\vert z \vert), j=1,...,n,
\end{eqnarray}
\begin{eqnarray}
\label{11'}
f_j(z) = o(\vert z \vert), j=n+1,...,N
\end{eqnarray}
and denote $'z' = (z'_1,....,z'_n)$, $''z' = (z'_{n+1},...,z'_N)$ so that $z' = ('z',''z')$.

For $w \in M \cap {\mathcal U}$ let

\begin{eqnarray*}
\sigma_w = \{ w' \in {\mathcal U}' \cap Q'_{f(w)}: \sum_{j=1}^N \rho'_j(f(w),\overline{w'}) T_kf_j(w,\overline w) =0, k=2,...,n \}
\end{eqnarray*}
It follows from (\ref{8}),(\ref{10}),(\ref{11}) that for $w \in M \cap {\mathcal U}$ the sets $\sigma_w$ are complex manifolds of dimension $N-n$ which smoothly depend
on $w$. Moreover, $f(w) \in \sigma_w$ and $T_{0'}\sigma_0 = \{ 'z' = 0 \}$. By (\ref{11}) we have $T_{0'}(S') = \{ ''z' = 0 \}$ and therefore $S'$ and $\sigma_0$ intersect transversally at $0'$. Thus $T_{f(w)}S'$ and $T_{f(w)}\sigma_w$ intersect also transversally and $S' \cap \sigma_{w} = \{ f(w) \}$. This implies (\ref{9}).

{\bf Remark.} Suppose  that the coordinates in $\C^N$ are ``normal'' for $M'$ at $0'$, i.e. the defining function of $M'$ can be chosen in the form

\begin{eqnarray}
\label{12}
\rho'(z',\overline z') = 2x_1' + \sum_{j=2}^N \vert z_j' \vert^2 + \sum_{\vert K \vert, \vert L \vert \geq 2} c_{KL}(y_1') \tilde z^K\overline{\tilde z^L}
\end{eqnarray}
where $\tilde z' = (z_2',...,z_N')$. Then $\sigma_0 = \{ 'z' = 0 \}$ and by lemma \ref{lemma2.1}

\begin{eqnarray}
\overline F \cap (\{ 0 \} \times {\mathcal U}' ) \subset \{ 0 \} \times \{ 'z'=0 \}
\end{eqnarray}

Set $\varphi_c(w,w') = \rho(w) + \rho'(w') - c[d(w',S')]^2$ and
$\Gamma_c = 
\{ (w,w') \in {\mathcal U} \times {\mathcal U}': \varphi_c(w,w') = 0 \}$. Since ${\mathcal U}$ 
and ${\mathcal U}$ are small, $\varphi_c \in C^\infty({\mathcal U} \times {\mathcal U}')$ and $\Gamma_c$ is a $C^\infty$ hypersurface passing through $(0,0')$. 

\begin{lemma}
\label{lemma2.3}
For any $c > 0$ the restriction of the Levi form of $\varphi_c$ to the complex tangent plane $T_{(0,0')}\Gamma_c$ has at least $2n-1$ positive eigenvalues.
\end{lemma}
\proof Since the tangent plane $T_{0'}(S')$ is an $n$-dimensional complex plane, the Levi form of the function $[d(w',S')]^2$ at $0'$ has $n$ zeros and $N-n$ positive 
eigenvalues. Thus for any $c >0$ the Levi form of $\varphi_c(w,w') = \rho(w) + \rho'(w') - c[d(w',S')]^2$ at $(0,0')$ has at least $n + N - (N-n) = 2n$ positive eigenvalues and its restriction to $T^c_{(0,0')}\Gamma_c$ has $\geq 2n-1$ positive eigenvalues.

Let $\Omega_c = \{ (w,w') \in {\mathcal U} \times {\mathcal U}': \varphi_c(w,w') > 0 \}$.

\begin{lemma}
\label{lemma2.4}
Let ${\mathcal U}$ and ${\mathcal U}'$ be small enough neighborhoods of $0$ and $0'$ respectively. For $c >0$ large enough the intersection $F \cap \Omega_c$ is 
closed in $\Omega_c$.
\end{lemma}

\proof If $(w,w') \in \overline F$ with $w \in M \cap {\mathcal U}$ then $w' \in Q'_{f(w)}$ and by corollary \ref{cor2} 
applied to $M'$ we obtain $\rho'(w') \leq c_1 \vert w' - f(w) \vert^2$. By lemma \ref{lemma2.2} $\vert w' - f(w) \vert^2 \leq c_2 [d(w',S')]^2$ 
and hence $\rho'(w') \leq c_1c_2[d(w',S')]^2$. Thus, if $(w,w') \in
{\mathcal U} \times {\mathcal U}'$ 
is a limit point for $F \cap \Omega_c$ and does not belong 
to $F$, then $\rho(w) = 0$, $\rho'(w') \leq c_1c_2[d(w',S')]^2$ and $(w,w')$ does not belong to $\Omega_c$ for $c \geq c_1c_2$.

\section{Reflection of analytic sets}

Let ${\mathcal U}$, ${\mathcal U}'$, $\rho$, $\rho'$, $M$, $M'$ be the
same as in the previous section and 
$(a,a') \in {\mathcal U} \times {\mathcal U}'$. We
can find an arbitrary small neighborhood $\Omega = \Omega(a,a') \subset {\mathcal U} \times {\mathcal U}'$ of $(a,a')$ and a neighborhood $V \times V' \subset {\mathcal U} \times 
{\mathcal U}'$ of $Q_a \times Q_{a'}$ such that for any $(w,w') \in V \times V'$ the intersection $(Q_w \times Q'_{w'}) \cap \Omega$ is non-empty and connected.

For such $\Omega$, a neighborhhod $V \times V'$ and a closed set $ A \subset \Omega$ we define its {\it reflection} $r(A)$ by

\begin{eqnarray}
\label{3.1}
r(A):= \{ (w,w') \in V \times V': S(w) \subset ({\mathcal U} \times Q'_{w'}) \cap \Omega \}
\end{eqnarray}
where $S(w):= (Q_w \times {\mathcal U}') \cap A$.

Notice that $r(A)$ depends not only on $A$ but also on $\Omega$ and $V \times
V'$. For  fixed $\Omega$, $V$ and $V'$, 
it follows immediately 
from (\ref{3.1}) that $\tilde A \subset A$ implies $r(A) \subset r(\tilde A)$.

If $(b,b') \in Q_a \times Q_{a'}$ and $\Omega(b,b')$ is an appropriate neighborhood of $(b,b')$ then we may consider the {\it second reflection} 
$r^2(A):= r(A_1)$ of $A_1:= r(A) \cap \Omega(b,b')$.

\begin{lemma}
\label{lemma3.1}
$A \subset r^2(A)$ near $(a,a')$.
\end{lemma}
\proof Let $(z,z') \in A$ be close enough to $(a,a')$. Then by (\ref{3.1}) it is enough to show that 
$A_1 \cap (Q_z \times {\mathcal U}') \subset  {\mathcal U} \times Q'_{z'}  $. Choose any point $(w,w') \in A_1 \cap (Q_z \times {\mathcal U}')$, i.e.
$(w,w') \in r(A) \cap \Omega(b,b')$ and $w \in Q_z$. Since $(z,z') \in A$ and $z \in Q_w$ it follows from (\ref{3.1}) that $z' \in Q'_{w'}$ and hence
$w' \in Q'_{z'}$. 

In this paper $A\subset \Omega$ will always be an analytic set. In general, its reflection $r(A)$ is not necessarily a (closed) analytic set in $V \times V'$. However analyticity of $r(A)$ can be established under certain additional conditions. In particular, if $\Omega = \omega \times \omega'$ and $A$ is the graph of a holomorphic map 
$g: \omega \longrightarrow \omega'$, then $r(A) \subset V \times V'$ is an analytic set defined by the condition $g(Q_w \cap \omega) \subset Q'_{w'}$. This and similar
cases have been previously discussed in  different papers (see, for instance, \cite{DiPi1,DiPi2}). The set $F$ introduced in section 2 of this paper 
is also a reflection of this kind.

\begin{lemma}
\label{lemma3.2}
Let $(a,a')$ be a point of an irreducible analytic set $A$ of dimension $d$ in a
neighborhood $\Omega = \Omega(a,a')$. 
Let $(b,b') \in Q_a \times Q_{a'}$ and $\dim S(b) = \dim 
(A \cap (Q_b \times {\mathcal U}')) = d-1$. Then there exists a neighborhood $\Omega(b,b')$ of $(b,b')$ such that after a possible shrinking of $\Omega(a,a')$ the set $r(A) \cap \Omega(b,b')$ is analytic
in $\Omega(b,b')$.
\end{lemma}
\proof Since $A$ is irreducible, the set $S(b) = A \cap (Q_b \times {\mathcal U}')$ is an analytic set in $\Omega(a,a')$ of pure dimension $d-1$. There exists a linear change of coordinates in $\C^n_z \times \C^N_{z'}$ such that in the new coordinates $(z^1,z^2) \in \C^{d-1} \times \C^{n+N-d+1}$ we have

\begin{eqnarray*}
S(b) \cap \{ z^1 = a^1 \} = \{ (a^1,a^2) \}
\end{eqnarray*}
where $(a^1,a^2)$ are the new coordinates of $(a,a')$. Consider a neighborhood $\Omega$ of $(a^1,a^2)$ of the form 
$\Omega = \Omega_1 \times \Omega_2 \subset \C^{d-1} \times \C^{n+N-d+1}$ such that $S(b)$ has no limit points on $\overline \Omega_1 \times \partial \Omega_2$. Then there exists a neighborhood $\Omega(b,b') = \omega(b) \times \omega'(b')$ such that $S(w)$ also does not have limit points on $\overline \Omega_1 \times \partial \Omega_2$ for any $w \in \omega(b)$ and therefore the projection $\pi:S(w) \longrightarrow \Omega_1$ is an $m$-sheeted branched holomorphic covering which depends antiholomorphically on $w$. There exists an open set $\omega_1 \subset \Omega_1$ such that $S(w) \cap (\omega_1 \times \Omega_2)$ is the union of the graphs of $m$ holomorphic mappings

\begin{eqnarray*}
z^2 = g^j(z^1,\overline w), z' = {g'}^j(z^1,\overline w), z^1 \in \omega_1, j=1,...,m
\end{eqnarray*}
These mappings also depend antiholomorphically on $w \in \omega(b)$.

By the uniqueness theorem the inclusion $S(w) \subset ({\mathcal U} \times Q'_{w'}) \cap \Omega$ is equivalent to the condition 
$S(w) \cap (\omega_1 \times \Omega_2) \subset ({\mathcal U} \times Q'_{w'}) \cap \Omega$ which can be expressed as 

\begin{eqnarray}
\label{3.2}
\rho'({g'}^j(z^1,\overline w),\overline w') = 0 
\end{eqnarray}
for all $z^1 \in \omega_1$ and $j=1,...,m$. This is a family of (anti)holomorphic equations for $w,w'$ and thus $r(A)$ is an analytic set in $\Omega(b,b')$.

\section{Proof of Theorem}

As in section 2 we assume that $\rho$, $\rho'$ and $f$ satisfy (\ref{10}),
(\ref{11}), (\ref{11'}). For 
any $w = (w_1,w_2,...,w_n) \in {\mathcal U}$ there exists unique ${}^sw = ({}^sw_1,{}^sw_2,...,{}^sw_n) \in Q_w$ such that $w_j = {}^sw_j$ for $j=2,...,n$. Since by \cite{Fors}
$f$ extends holomorphically across a dense open set $\Sigma \subset M$, 
it is holomorphic on some open set ${\mathcal U}^-_1$ containing ${\mathcal U}^- \cup \Sigma$. There also exists 
an open set ${\mathcal U}^+_1$ containing ${\mathcal U}^+ \cup \Sigma$ such that ${}^sw \in {\mathcal U}^-_1$ for any $w \in {\mathcal U}^+_1$. Denote by $Q^c_w$ the connected component of $Q_w \cap {\mathcal U}^-_1$ that contains ${}^sw$. We can now modify the definition of $F$ and consider

\begin{eqnarray}
\label{4.1}
 F_1:= \{ (w,w') \in {\mathcal U}^+_1 \times {\mathcal U}': f(Q^c_w) \subset Q'_{w'} \}
\end{eqnarray}
Obviously $ F_1$ coincides with $F$ over ${\mathcal U}^+$. The proof of lemma \ref{lemma1.4} works for $ F_1$ without any changes and thus 
$ F_1$ is an analytic set in ${\mathcal U}^+_1 \times {\mathcal
  U}'$. The intersection ${\mathcal U}_1: = 
{\mathcal U}_1^+ \cap {\mathcal U}^-_1$ is an open 
neighborhood of $\Sigma$ and $ F_1$ contains the graph of $f$ over ${\mathcal U}_1$. 

The set $ F_1$ consists of irreducible components of two types. We say that a component of $ F_1$ is {\it relevent} 
if it contains an open piece of the graph of $f$ over ${\mathcal U}_1$. Otherwise we call it {\it irrelevent}. Thus 
$ F_1$ is the union of two analytic sets: $F_r$ and $F_i$ which consist of all  relevent and  irrelevent
components respectively.  It is obvious that the dimension of $F_r$ is $\geq n$ at any its point, the dimension of the intersection of $F_i$ with the graph of $f$ is $<n$ and $(0,0') \in \overline F_r$.

We now represent $F_r$ as $F_r^{(n)} \cup F_r^{(n+1)}$ where $F^{(n)}_r$ is the union of all $n$-dimensional relevent components and $F_r^{(n+1)}$ consists of all relevent components of dimension $\geq n +1$.

There are two possibilities:
\begin{itemize}
\item[(1)] After shrinking ${\mathcal U}$ and ${\mathcal U}'$ we have $F_r = F_r^{(n)}$. 
\item[(2)] $(0,0') \in \overline{F_r^{(n+1)}}$.
\end{itemize}

We first prove Theorem in the second case.

\subsection{Proof of Theorem in the case (2)} 
We need the following technical statement which is a slight  variation
of the standard results (see, for instance, \cite{Ch}, p. 36).
\begin{lemma}
\label{convergence}
Let $A$ be a complex purely $m$-dimensional analytic set in a domain $\Omega
\subset \C^n$ and $(A_\nu)$ be a sequence of purely $p$-dimensional 
complex analytic sets in $\Omega$. Suppose that $p \geq m$ and that
the cluster set $cl(A_\nu)$ of the sequence $(A_\nu)$ is contained in
$A$. 
Then $ p = m$ and $cl(A_\nu)$ is a  union of some irreducible
components of $A$.
\end{lemma}
As usual, by the cluster set $cl(A_\nu)$ of a sequence $(A_\nu)$ we mean the
set of all points $a \in \Omega$ such that there exists a subsequence
$(\nu(k))$ of
indices and points
$a_{\nu(k)} \in A_{\nu(k)}$ converging to $a$ as $k$ tends to infinity.  
For the convenience of readers we give the proof of the lemma.
\proof Fix a point $a \in cl(A_\nu)$; we can assume that $a = 0$. 
Consider a complex linear $(n-m)$-dimensional subspace $L$ of $\C^n$
satisfying   $A \cap L = \{ 0 \}$. Then
there exist a ball $B$ centered at the origin and $r > 0$ such that 
the distance form $A$ to $L \cap \partial B$ is equal to $r$. Since 
$cl(A_\nu) \subset A$, for every $\nu$ big enough the sets $A_\nu$
do not intersect the $r/2$-neighborhood of $L \cap \partial B$.
On the other hand, $0 = \lim a_{\nu(k)}$ with $a_{\nu(k)} \in A_{\nu(k)}$  so for any $k$ big enough the
intersection $A_{\nu(k)} \cap B$ is not empty.   Then the intersection $(L +
a_{\nu(k)}) \cap A_{\nu(k)} \cap B$ is a compact analytic subset in $B$ and so
its dimension is equal to $0$. Since $\dim L = n-m$, this implies $p = \dim
A_{\nu(k)} \leq m$ and we obtain that $m = p$. 

Now we  prove that $cl(A_\nu)$ coincides with a union of some irreducible
components of $A$. Since the set ${\mathcal S}(A)$ of singular points of $A$
is  an analytic set of dimension $< m$, it follows from the first part of
lemma that $cl(A_\nu)$ is not contained in ${\mathcal S}(A)$. So the
intersection
of $cl(A_\nu)$ with the set ${\mathcal R}(A)$ of regular points of $A$ is not
empty and  this is sufficient to
show that $cl(A_\nu)$ is open in ${\mathcal R}(A)$. Consider an arbitrary
point    $a \in cl(A_\nu) \cap {\mathcal R}(A)$.
As above, we assume that $a = 0$. After a biholomorphic change of
coordinates we can assume that in a neighborhood of the origin $A$
coincides with the coordinate space $P$ of variables $z_1,...,z_p$.
Denote by $L$ the coordinates space of variables $z_{p+1},...,z_n$
and fix small enough the balls $B \subset P$ and $B'\subset L$
centered at the origin. Since $cl(A_\nu) \subset A$, for every $\nu$
big enough the set $A_\nu$ does not intersect $B \times \partial B'$.
So every $A_\nu \cap (B \times B')$ is a  analytic covering brunched
over $B$. Hence for every point $b \in B$ the fiber $\{ b \} \times L$ 
 contains a point $(b,c_\nu) \in A_\nu \cap (B \times B')$.
Since $cl(A_\nu) \subset A$, we get $\lim c_\nu = 0$  which
proves the claim.

\begin{lemma}
\label{lemma4.1}
If $(0,0') \in \overline{F_r^{(n+1)}}$, then $F_r^{(n+1)}$ extends to an analytic set in a neighborhood of $(0,0')$.
\end{lemma}
\proof Since $F_r^{(n+1)}$ contains only the relevant  components, 
there exists a sequence $w^\nu \in \Sigma$ converging to $0$ as $\nu
\longrightarrow \infty$ such that $(w^\nu,f(w^\nu)) \in F_r^{(n+1)}$ for any
$\nu$ (if not, the proof is reduced to the case (1)). 
Let $\tilde f$ be a $C^\infty$ extension of $f$ to ${\mathcal U}$ that
coincides with $f$ on ${\mathcal U}^-_1$ and 
$S' = \tilde f({\mathcal U}) \subset {\mathcal U}'$. Let $\varphi_c$,
$\Gamma_c$ and $\Omega_c$ be the same as in section 2. Since $d(w',S')
= 0$ for $w' = f(w)$, 
$w \in {\mathcal U}^-_1 \cap {\mathcal U}^+$ the intersection 
$F_r^{(n+1)} \cap \Omega_c$ is not empty and moreover $(0,0') \in
\overline{F_c^{(n+1)} \cap \Omega}$. 
The set $F_r^{(n+1)}$ can be decomposed to a finite union of
analytic sets (perhaps, reducible) of pure dimensions $\geq n +1$:
\begin{eqnarray*}
F_r^{(n+1)} = \cup_{k \geq n+1} F_{r,k}^{(n+1)}, \dim F_{r,k}^{(n+1)} = k
\end{eqnarray*}
(see, for instance, \cite{Ch},
p.51). By lemma \ref{lemma2.3} the Levi form of $\varphi_c$ has at least $2n-1$
positive eigenvalues on $T^c_{(0,0')}\Gamma_c$. 
Since $N < 2n$, the set $F_r^{(n+1)}$ extends to an analytic set in a neighborhood of
$(0,0')$ by Rothstein's theorem on the analytic extension across
pseudoconcave hypersurfaces (see, 
for instance, \cite{Ch}). More precisely,   there exists an analytic set $ 
 \tilde F \subset {\mathcal U} \times {\mathcal U}'$ such that 
$ F_r^{(n+1)} \subset \tilde F \cap ({\mathcal U}_1^+ \times {\mathcal U}')$. 
This set contains the graph of $f$ near $(0,0')$.   Every irreducible  component of $\tilde
F$ is of the dimension   $\geq n +1$ and has  a non-empty open subset  contained in $F_r^{(n+1)}$.

\begin{proposition}
\label{pro4.2}
If $(0,0') \in \overline{F_r^{(n+1)}}$, then $f$ extends holomorphically to a neighborhood of $0$.
\end{proposition}

 We begin the proof with the following 

\begin{lemma}
\label{lemma4.3}
In any neighborhood of $(0,0')$ there exists a point $(w^0,{w'}^0) \in F_r^{(n+1)} \cap (Q_0 \times Q_{0'}')$ with $w^0 \neq 0$, ${w'}^0 \neq 0$. Moreover $(0,0')$ belongs to the closure 
of some component of $F_r^{(n+1)} \cap (Q_0 \times Q_{0'}')$ which contains $(w^0,{w'}^0)$.
\end{lemma}
\proof Suppose that the coordinates $z' = ('z',''z')$ in $\C^N$ are ''normal'' for $M'$ at $0'$, i.e. satisfy (\ref{12}).
Consider a sequence $(a^\nu, {a'}^\nu) \in F_r^{(n+1)}$ such that
$(a^\nu,{a'}^\nu) \longrightarrow (0,0')$ and for every
$\nu$ we have  $a^\nu \in
{\mathcal U}_1^- \cap {\mathcal U}^+$, ${a'}^\nu = f(a^\nu)$. Passing to a
subsequence, we may also assume that there exists an irreducible component of
$\tilde F$ of dimension $d \geq n + 1$ containing the graph of $f$ in a neighborhood
of $(0,0')$ such that  $(a^\nu,{a'}^\nu)$ belongs to this component for every
$\nu$. We denote it again by $\tilde F$.  Let $b^\nu = {}^sa^\nu$. For any $\nu = 1,2,...$ the intersection 
$\tilde S_\nu:= \tilde F \cap (Q_{b^\nu} \times {\mathcal U}')$ is an
analytic set in ${\mathcal U} \times {\mathcal U}'$ of pure dimension $d-1$
and containing $(a^\nu,{a'}^\nu)$. Indeed, if not, $\tilde F$ is contained
in the hypersurface $Q_{b^\nu} \times {\mathcal U}'$ which is impossible since
$\tilde F$ contains an open piece of the graph of $f$. For the same reason
the dimension of the set $\tilde F \cap (Q_0 \times {\mathcal U}')$
also is equal to $d - 1$. 
The cluster set $\tilde S_0:= cl(\tilde S_\nu)$ of the sequence $\tilde S_\nu$
with respect to ${\mathcal U} \times {\mathcal U}'$ is contained in $\tilde F
\cap (Q_0 \times
{\mathcal U}')$ and so by lemma \ref{convergence} 
 $\tilde S_0$ is the union of some components of 
$\tilde F \cap (Q_0 \times {\mathcal U}')$; in particular, $\tilde S_0$ is an
analytic set of pure dimension $d-1$. 
On the other hand, denote by $F_r^d$ the union of purely $d$-dimensional
components of $F_r^{(n+1)}$  having a non-empty open intersection with $\tilde F$.  
For any $\nu = 1,2,...$ the intersection 
$S_\nu:= F_r^{d} \cap (Q_{b^\nu} \times {\mathcal U}')$ is an analytic set
in ${\mathcal U}^+ \times {\mathcal U}'$ of pure dimension $d-1$ 
containing the point $(a^\nu,{a'}^\nu)$ and contained in $\tilde S_\nu$. Since
$Q_{b^\nu} \subset {\mathcal U}_1^+$  and $F_r^{(n+1)}$ is an analytic set in ${\mathcal U}_1^+ \times
{\mathcal U}'$, every $S_\nu$ is a (closed) analytic subset in ${\mathcal U}
\times {\mathcal U}'$ and so coincides with a union of some irreducible
components of $\tilde S_\nu$. Therefore by lemma \ref{convergence} the cluster set 
$S_0:= cl(S_\nu) \subset \overline{F_r^{(n+1)}} \cap (Q_0 \times {\mathcal
  U}')$  
is the union of some components of 
$\tilde S_0$ and  $\dim S_0 = d-1 \geq n$. 
 
On the other hand  by the remark after lemma \ref{lemma2.2} 
$\overline{F_r^{(n+1)}} \cap ( \{ 0 \} \times {\mathcal U}') \subset
\{ 0 \} \times \sigma_0$, where 
$\sigma_0 = \{ z' \in \C^N: 'z'=0 \}$. Since $\dim S_0 = d-1 \geq n$ and 
$\dim \sigma_0 = N - n < n$, the set $S_0$ is not contained in $\{ 0 \} \times
\sigma_0$. Therefore, 
 $S_0$ is not contained in  $ \overline{F_r^{(n+1)}} \cap ( \{ 0 \} \times
{\mathcal U}')$.
Hence  in any neighborhood of $(0,0')$  there exists a point $ (w^0,{w'}^0) \in
S_0$ with $w^0 \neq 0$. Moreover, the set $S_0$ is not contained in $Q_0
\times \{ 0' \}$   because $\dim S_0  > n-1 = \dim
Q_0$. Therefore, in any neighborhood of $(0,0')$ there exists a point  
$ (w^0,{w'}^0) \in
S_0$ with $w^0 \neq 0$, ${w'}^0 \neq 0$. 
Moreover $w^0 \in {\mathcal U}^+$ because $Q_0 \subset {\mathcal U}^+ \cup \{ 0 \}$. This means
that $(w^0,{w'}^0) \in {\mathcal U}^+ \times {\mathcal U}'$ and thus
$(w^0,{w'}^0) \in {F_r^{(n+1)}} \cap S_0$. Finally, it follows from the definition
(\ref{4.1}) that every $S_\nu$ is contained in ${\mathcal U} \times
Q'_{f(b^\nu)}$. Hence, $S_0 \subset Q_0 \times Q'_{0'}$ and we get the first
claim of lemma.

Prove the second claim. After possible shrinking of ${\mathcal U}$ and
${\mathcal U}'$ the set $\tilde F \cap (Q_0 \times Q'_{0'}) \times ({\mathcal
  U} \times {\mathcal U}')$ consists of  a finite number of ireducible
components and every such  component contains $(0,0')$. So the second claim
follows from the first part of lemma.

{\it Proof of proposition \ref{pro4.2}:} Consider a sequence of points
$(w^\nu,{w'}^\nu) \in F_r^{(n+1)} \cap ( Q_0 \times Q'_{0'})$ such that
$(w^\nu, {w}'^\nu) \longrightarrow (0,0')$ and $w^\nu \neq 0$, ${w'}^\nu \neq
0$. Choose appropriate neighborhoods $\Omega_\nu$ and $\Omega_\nu^0$ of $(w^\nu,{w'}^\nu)$ and $(0,0')$ respectively such that $F^2_\nu:=
r(F_r^{(n+1)} \cap \Omega_\nu)$ is an analytic set in $\Omega_\nu^0$. Consider
the analytic sets $A_\mu: = \cap_{\nu = 1}^\mu F^2_\nu$. Then $A_{\mu +1} \subset
A_\mu$. If $d_\mu$ denotes the dimension of $A_\mu$, then $d_{\mu + 1} \leq
d_{\mu}$.   Therefore $d_\mu \geq
n$ for every $\mu$ and there exists $d_{\mu_0} \geq n$ such that $d_\mu =
d_{\mu_0}$ for any $
\mu > \mu_0$.  Since $F_1$ defined by (\ref{4.1}) is the reflection of $\Gamma_f$
every set $F^2_\nu$ contains an open 
piece of the graph $\Gamma_f$ by lemma \ref{lemma3.1}.
 Since the set $A_{\mu_0}$ has a finite number of components,
there exists a neighborhood $\Omega^0$ of $(0,0')$ and an irreducible analytic
set $F^2 \subset \Omega^0$ containing an open piece of $\Gamma_f$ such that
$\dim F^2 = d_{\mu_0} \geq n$ and for any $\nu$ big enough $F^2 \cap
\Omega^0_\nu \subset F^2_\nu$.  
Since $F_2$ is an
analytic set in a neighborhood of $(0,0')$, we can consider its reflection
$F_3:= r(F_2)$ which is an analytic set in a neighborhood $\tilde \Omega^0$ of
$(0,0')$. For any $\nu$ big enough $(w^\nu,{w'}^\nu) \in \tilde \Omega^0$.
 By lemma \ref{lemma3.1} $F_3$  contains $F_r^{(n+1)}$ near
 $(w^\nu,{w'}^\nu)$. Recall that $F_r^{(n+1)}$ has a finite number of
 irreducible components in a neighborhood of the origin and every component
contains an open
 piece of $\Gamma_f$.
Hence, $F_3$ contains an open piece of $\Gamma_f$ in a neighborhood of the origin.
  Thus both $F_2$ and $F_3$ are analytic sets in a neighborhood of $(0,0')$ and both contain a piece of $\Gamma_f$ as well as $F_2 \cap F_3$.

\begin{lemma}  
\label{claim}
We have $F_2 \cap F_3 \cap ( \{ 0 \} \times {\mathcal U}') = \{ (0,0') \}$. 
\end{lemma}
Indeed, suppose that $\gamma:= F_2 \cap F_3 \cap ( \{ 0 \} \times
{\mathcal U}')$ is an analytic set 
of positive dimension. Let $(0,{z'}^0) \in F_2 \cap F_3$. We have 
$F_3 = \{ (z,z'): F_2 \cap (Q_z \times {\mathcal U}') \subset
{\mathcal U} \times Q_{z'}' \}$. Hence 
$F_2 \cap (Q_0 \times {\mathcal U}') \subset {\mathcal U} \times
Q_{{z'}^0}'$ and thus ${z'}^0 \in Q_{{z'}^0}$ that 
is ${z'}^0 \in M'$. Hence $\gamma \subset M'$ which contradicts  the strict pseudoconvexity of $M'$. This proves the claim.

Therefore $F_2 \cap F_3$ has a locally proper at $(0,0')$ projection
$\pi: F_2 \cap F_3 \longrightarrow {\mathcal U}$. 
Hence $\dim (F_2 \cap F_3) = n$ and $F_2 \cap F_3$ 
is an analytic
continuation of $\Gamma_f$ to a neighborhood of $(0,0')$. By \cite{BB}
$f$ extends holomorphically to a neighborhood of $(0,0')$. This
completes the proof of proposition \ref{pro4.2} and proves theorem in the case (2).

\subsection{Proof of Theorem in the case (1)} Consider now the case (1) where
$\dim F_r = n$ is a neighborhood of the origin. Everywhere below we suppose
that this assumption holds.

\begin{lemma} Suppose that $$\overline F_r^{(n)} \cap (Q_0 \times Q'_{0'}) = 
\overline F_r^{(n)} \cap (\{ 0 \} \times Q'_{0'}).$$ 
Then $$\overline F_1 \cap (\{ 0
\} \times Q'_{0'}) = \{
  0 \} \times \sigma_0.$$
\end{lemma} 
\proof  Let $\overline F_r^{(n)} \cap (\{ 0 \} \times Q'_{0'}) = \{
  0 \} \times X$. It follows by lemma \ref{lemma2.1} that $X \subset \sigma_0$.
 Consider a sequence $(w^\nu)$ of points 
in $\Sigma$ converging to $0$ and set ${w'}^\nu = f(w^\nu) \in
M'$. Denote by ${}_{w^\nu}Q_z$ the germ of the Segre variety $Q_z$ at
$w^\nu$ and consider the analytic  sets 
\begin{eqnarray*}
S_\nu = \{ (z,z') \in Q_{w^\nu} \times Q_{{w'}^\nu} : f({}_{w_{\nu}}Q_z) \subset Q'_{z'} \}
\end{eqnarray*}
in ${\mathcal U} \times {\mathcal U}'$. Then $\dim S_\nu \geq n - 1$.
Since 
$ \dim \sigma_0 = N - n \leq 2n - 1 -n = n-1$ and $cl(S_\nu) \subset
\{ 0 \} \times X$ (by the hypothesis of lemma) , lemma \ref{convergence} implies that 
$$cl(S_\nu) = \{ 0 \} \times X = \{ 0 \} \times \sigma_0$$
 which  proves lemma. 

 We claim that the projection of $\overline F_r^{(n)} \cap (Q_0 \times
Q'_{0'})$ to $\C^n$ can not be equal to the singleton $\{ 0 \}$. Indeed,
assume that  $\overline F_r^{(n)} \cap (Q_0 
\times Q'_{0'}) = \{ 0 \} \times \sigma_0$. Fix  a point ${z'}^0 \in \sigma_0$
which does not belong to $M'$ (since $M'$ is
strictly pseudoconvex, it contains no analytic sets of positive
 dimension). Consider a sequence of points $(z^\nu,{z'}^\nu) \in
 S_\nu$ converging to $(0,{z'}^0)$. Consider analytic sets 
$A_\nu =  F_r^{(n)} \cap (Q_{z^\nu} \times Q_{{z'}^\nu})$.
Since $F_r^{(n)}$ contains the graph of
 $f$ over $\Sigma$, for every $\nu$ we have  $\dim A_\nu \geq n -1$. 
We have $cl(A_\nu) \subset \{ 0 \} \times \sigma_0$. Hence lemma
\ref{convergence} implies $cl(A_\nu) = \{ 0 \} \times \sigma_0$. On
the other hand, $cl(A_\nu) \subset \{ 0 \} \times Q'_{{z'}^0}$. 
 Hence ${z'}^0 \in \sigma_0 \subset Q'_{{z'}^0}$ and so ${z'}^0 \in M'$:
 a contradiction. 
 
Thus, in any neighborhood of $(0,0')$ the intersection $\overline F_r^{(n)} \cap
(Q_0 \times Q'_{0'})$ 
contains points $(z,z')$ with $z \neq 0$. Let us show that for every such
point we also have $z' \neq 0$. Indeed, assume by contradiction that 
$(z,0) \in \overline F_r^{(n)} \cap (Q_0 \times Q'_{0'})$ and $z \neq
0$. Therefore $z \in {\mathcal U}^+$  and $(z,0) \in F_r^{(n)}$. Then 
$f(Q_z \cap {\mathcal U}^-) \subset Q'_{0'} \subset {{\mathcal U}'}^+ \cup \{
0' \}$. On the other hand,  
$f(Q_z \cap {\mathcal U}^-) \subset {{\mathcal U}'}^-$ and $Q'_{0'} 
\cap {{\mathcal U}'}^- = \{ 0 \}$. This implies that $f$ vanishes
identically on the complex hypersurface $Q_z \in {\mathcal U}^-$ (we point out
that $Q_z$ intersects $M$ transversally at the origin since $z \in Q_0$ and $z
\neq 0$). However, $f$ has the maximal rank: a contradiction.  

We sum up this considerations in the
following statement.

\begin{lemma}
\label{nondegeneracy}
In any neighborhood of $(0,0')$ there exists a point $(w^0,{w'}^0) \in 
\overline F_r^{(n)} \cap (Q_0 \times Q'_{0'})$ with $w^0 \neq 0$ and 
${w'}^0 \neq 0$.
\end{lemma}

Now we are able to conclude the proof of Theorem. Fix 
 a point $(w^0,{w'}^0) \in F_r^{(n)} \cap (Q_0 \times Q'_{0'})$ with $w^0 \neq 0$ and 
${w'}^0 \neq 0$, fix a neighborhood $\Omega_0$ of $(w^0,{w'}^0)$ and consider the
reflection $F_2 = r(F_r^{(n)} \cap \Omega_0)$. Then $F_2$ is an analytic set in a
neighborhood $\Omega_1$ of $(0,0')$  and contains 
an open piece of $\Gamma_f$; in particular, $\dim F_2 \geq n$. If $\dim F_2 =
n$, we conclude. If not, we apply  an argument similar to 
 the proof of the case (2).

Consider   a basis ${\mathcal U}_\nu \times {\mathcal U}'_\nu$ of
neighborhoods of $(0,0')$ and a sequence $(w^\nu,{w'}^\nu) \in  \overline
F_r^{(n)} \cap (Q_0 \times Q'_{0'}) \cap ({\mathcal U}_\nu \times {\mathcal U}'_\nu)$ with $w^\nu \neq 0$ and 
${w'}^\nu \neq 0'$ such that for every $\nu$ there exists a component of 
$F_r^{(n)} \cap ({\mathcal U}_\nu \times {\mathcal U}'_\nu)$ containing the
point $(w^\nu,{w'}^\nu)$ and an open subset of $\Gamma_f$ in ${\mathcal U}_\nu \times {\mathcal U}'_\nu$.
 Choose appropriate neighborhoods $\Omega_\nu$ and $\Omega_\nu^0$ of $(w^\nu,{w'}^\nu)$ and $(0,0')$ respectively such that $F^2_\nu:=
r(F_r^{(n)} \cap \Omega_\nu)$ is an analytic set in $\Omega_\nu^0$. As in the
proof of the case (2), consider
the analytic sets $A_\mu: = \cap_{\nu = 1}^\mu F^2_\nu$.  As above let $d_\mu$ denotes the dimension of $A_\mu$. Then $d_{\mu + 1} \leq
d_{\mu}$ and $d_\mu \geq
n$ for every $\mu$; furthermore, there exists $d_{\mu_0} \geq n$ such that $d_\mu =
d_{\mu_0}$ for any $
\mu > \mu_0$.  Since the set $F_1$ defined by (\ref{4.1}) is the reflection of
graph of $f$
every set $F^2_\nu$ contains an open 
piece of  $\Gamma_f$ in view of lemma \ref{lemma3.1}.
 The set $A_{\mu_0}$ has a finite number of components, so
there exists a neighborhood $\Omega^0$ of $(0,0')$ and an irreducible analytic
set $F^2 \subset \Omega^0$ containing an open piece of $\Gamma_f$ such that
$\dim F^2 = d_{\mu_0} \geq n$ and for any $\nu$ big enough $F^2 \cap
\Omega^0_\nu \subset F^2_\nu$.  
Since $F_2$ is an
analytic set in a neighborhood of $(0,0')$, we can consider its reflection
$F_3:= r(F_2)$ which is an analytic set in a neighborhood $\tilde \Omega^0$ of
$(0,0')$. For any $\nu$ big enough $(w^\nu,{w'}^\nu) \in {\mathcal U}_\nu
\times {\mathcal U}'_\nu \subset \tilde \Omega^0$.
 By lemma \ref{lemma3.1} $F_3$  contains $F_r^{(n)}$ near
 $(w^\nu,{w'}^\nu)$. Hence $F_3$ contains any component of $F_r^{(n)} \cap
 ({\mathcal U}_\nu \times {\mathcal U}'_\nu)$ 
passing through $(w^\nu,{w'}^\nu)$.
Therefore, $F_3$ contains an open piece of $\Gamma_f$.
  We obtain that both $F_2$ and $F_3$ are analytic sets in a neighborhood of $(0,0')$
  and both contain an open piece of $\Gamma_f$ as well as $F_2 \cap F_3$.
   Repeating the proof of lemma
\ref{claim} we get that the projection  $\pi: \tilde F_2 \cap F_3 \longrightarrow {\mathcal U}$
is locally proper at $(0,0')$ and we conclude as in the case (2).

\bigskip 
\bigskip

{\footnotesize

Sergey Pinchuk: Department of Mathematics, Indiana University,  Bloomington
IN-47405, USA and South Ural State University ,76 Lenin ave., Chelyabinsk, 45080 , Russia, 
  {\rm pinchuk@indiana.edu} 

Alexandre Sukhov: Universite des Sciences et Technologies de Lille, U.F.R. de
Mathematiques, Laboratoire Paul Painleve, 59655, Villeneuve d'Ascq, Cedex,
France, {\rm sukhov@math.univ-lille1.fr}

}

\end{document}